\documentclass[a4paper,12pt]{amsart}

\pdfoutput=1

\usepackage[foot]{amsaddr}
\usepackage{a4wide}
\usepackage[T1]{fontenc}
\usepackage[utf8]{inputenc}

\usepackage{amsmath, amsfonts, amsthm}
\usepackage{mathtools}
\usepackage[svgnames]{xcolor}

\definecolor[named]{LightGray}{rgb}{0.95,0.95,0.96}

\usepackage{listings}

\definecolor{codegreen}{rgb}{0,0.6,0}
\definecolor{codegray}{rgb}{0.5,0.5,0.5}
\definecolor{sageblue}{RGB}{94, 94, 255}
\definecolor{codepurple}{rgb}{0.58,0,0.82}
\definecolor{backcolour}{rgb}{0.98,0.98,0.96}

\lstdefinestyle{mystyle}{
    frame=tb,
    language=python,
    keywordstyle=\color{magenta},
    escapeinside={@}{@},
    backgroundcolor=\color{backcolour},   
    commentstyle=\color{codegreen},
    numberstyle=\tiny\color{codegray},
    stringstyle=\color{codegreen},
    basicstyle=\ttfamily\footnotesize,
    breakatwhitespace=false,         
    breaklines=true,                 
    captionpos=b,                    
    keepspaces=true,                 
    numbersep=5pt,                  
    showspaces=false,                
    showstringspaces=false,
    showtabs=false,                  
    tabsize=2,
    upquote=true,
}
\usepackage{hyperref}

\DeclarePairedDelimiter{\abs}{\lvert}{\rvert}

\newif\ifdetails
\detailstrue
\newcommand{\DETAIL}[1]%
{\ifdetails\par\fbox{\begin{minipage}{0.9\linewidth}\textit{Detail:}
      #1\end{minipage}}\par\fi}
\newcommand{\TODO}[1]%
{\ifdetails\par\fbox{\begin{minipage}{0.9\linewidth}\textbf{TODO:}
      #1\end{minipage}}\par\fi}

\newtheorem{lemma}{Lemma}

\newtheorem{theorem}[lemma]{Theorem}
\theoremstyle{remark}

\allowdisplaybreaks

\DeclareMathOperator{\res}{Res}
\newtheorem{problem}{Problem}

\newcommand{\code}[1]{\texttt{\detokenize{#1}}}

\begin{document}
\title[Bivariate asymptotics with explicit error terms]{A computer algebra package for bivariate asymptotics with explicit error terms}

\author{Benjamin {Hackl}}
\address{Department of Mathematics and Scientific Computing, University of Graz, Heinrichstraße 36, 8010 Graz, Austria}
\email{benjamin.hackl@uni-graz.at}
\thanks{
  This article is the full version of an extended abstract~\cite{Hackl-Wagner:2024:binomial-mellin-AofA24} that was presented
  at the 35th International Conference on Probabilistic, Combinatorial
  and Asymptotic Methods for the Analysis of Algorithms (AofA 2024).
}

\author{Stephan {Wagner}}
\address{Institute of Discrete Mathematics, TU Graz, Steyrergasse 30, 8010 Graz, Austria \and Department of Mathematics, Uppsala University, Box 480, 751 06 Uppsala, Sweden}
\email{stephan.wagner@tugraz.at}
\thanks{Stephan Wagner is supported by the Swedish research council (VR), grant 2022-04030.}


\begin{abstract}
Making use of a newly developed package in the computer mathematics system SageMath, we show how to perform a full asymptotic analysis of certain types of sums that occur frequently in combinatorics, including explicit error bounds. We present two applications of the general approach to illustrate its use: the first concerns a classical problem due to Ramanujan, while the second one concerns a question of B\'ona and DeJonge on 132-avoiding permutations with a unique longest increasing subsequence that can be translated into an inequality for a certain binomial sum.
\end{abstract}

\maketitle

\section{Introduction}
\label{sec:intro}

Sums involving binomial coefficients and other hypergeometric expressions occur frequently in enumerative and analytic combinatorics. For example,
\begin{equation}\label{eq:example1}
\sum_{k = 0}^n \frac{1}{k+1} \binom{n+k}{n} \binom{n}{k}
\end{equation}
yields the large Schr\"oder numbers, which count (among other things) many different types of lattice paths and permutations. The sum
\begin{equation}\label{eq:example2}
\sum_{k=0}^{\lfloor n/2 \rfloor} \binom{n}{2k} \frac{(2k)!}{2^k k!}
\end{equation}
counts involutions, or matchings in complete graphs.
There is a well-established mathematical toolkit for dealing with such sums,
based on techniques such as the \emph{(discrete) Laplace method},
the \emph{Stirling approximation} of factorials and binomial coefficients,
and the \emph{Mellin transform}. See~\cite{Flajolet2009analytic} for a comprehensive
account of these and many other tools.

While these methods are well known and in some sense mechanical, it is still not
straightforward to implement them in a computer as they often involve
ad-hoc estimates and careful splitting into different cases/regions of summation
that are analyzed separately. Moreover, while a lot of the complications can
often be hidden in $O$-terms, things become more involved when explicit error
bounds are desired.

This paper aims to make a contribution towards building a toolkit for asymptotic
analysis of the type of sums exemplified by~\eqref{eq:example1} and~\eqref{eq:example2}. This is done by means of a computer algebra package capable of performing bivariate asymptotics,
including automatically derived, guaranteed error bounds with explicit constants. This new package, \code{dependent_bterms},
for the computer mathematics system
SageMath~\cite{SageMath} enhances the core implementation of asymptotic expansions
and in particular the arithmetic with SageMath's analogue of $O$-terms with explicit
error bounds, called $B$-terms. The package
takes care of the ``heavy lifting'' when it comes
to the mechanical computations of asymptotic expansions with explicit error bounds,
and allows (in contrast to the functionality already included in SageMath)
computations with an additional variable $k = k(n)$ that depends on the ``main''
variable $n$. An in-depth introduction to the package can be found in
Section~\ref{sec:bterms}.

After introducing the package, we choose two illustrative examples to demonstrate how
tedious mechanical computations can be delegated to our package. In
Section~\ref{sec:ramanujan-example} we study a classical example in form
of the Ramanujan $Q$-function. This function and its ``complement'' are defined by
\[
  Q(n) = \sum_{k=0}^{n-1} \frac{n!}{(n-k-1)!\, n^{k+1}}
  = 1 + \frac{n-1}{n} + \frac{(n-1)(n-2)}{n^2} + \cdots + \frac{(n-1)!}{n^{n-1}},
\]
and
\[
  R(n) = \sum_{k=0}^{\infty} \frac{n!\, n^k}{(n+k)!}
  = 1 + \frac{n}{n+1} + \frac{n^2}{(n+1)(n+2)} + \cdots.
\]
The function $R(n)$ can be seen as a complement of $Q(n)$ in the sense that
it is not hard to verify the relation
\[ Q(n) + R(n) = \frac{n! e^n}{n^n} \]
by splitting the power series for $e^n$ into two parts. One of the problems posed
by Ramanujan can be formulated in terms of $Q(n)$ and $R(n)$ as follows.

\begin{problem}[Ramanujan, Question 294~\cite{Ramanujan:1911:question-294}]\label{prob:ramanujan}
  Show that the difference between $R(n)$ and $Q(n)$ can be expressed as
  \begin{equation*}
    R(n) - Q(n) = \frac{2}{3} + \frac{8}{135(n + \vartheta)}
  \end{equation*}
  for some $\vartheta = \vartheta(n)$ that satisfies $2/21 \leq \vartheta(n) \leq 8/45$.
\end{problem}

This is not actually an open problem, as it was settled completely by Flajolet, Grabner,
Kirschenhofer and Prodinger~\cite{Flajolet-Grabner-Kirschenhofer-Prodinger:1995:ramanujan-q},
following earlier work of Ramanujan, Watson and
Knuth~\cite{Knuth:1973:TAOCPv3, Watson:1929:theorems-ramanujan}. It is, however,
a nice and classical example for a question whose answer involves precise asymptotic
expansions with quantified error terms. Moreover, the high precision in Ramanujan's statement makes it
particularly challenging.

As already mentioned above,
sums like the one underlying this problem are well-studied with different
established methods that allow extracting the asymptotic growth. In this article,
we focus on the Mellin transform framework (see
\cite[Section 5]{Flajolet-Gourdon-Dumas:1995:mellin} for a rigorous treatment),
which roughly proceeds along the following lines:
\begin{itemize}
\item Split the sum into ``small'' and ``large'' regions, depending on the value of the
  summation index $k$.
\item Show that the contribution of large values is negligible.
\item Approximate the summands in a suitable way and ensure that the (sum over the) approximation error
  is negligibly small.
\item Turn the sum into an infinite sum, again at the expense of a negligible error term.
\item Apply the Mellin transform to obtain an integral representation for the resulting infinite sum.
\item Use residue calculus to derive the final asymptotic formula.
\end{itemize}

As we will see, carrying out these steps and using the \code{dependent_bterms} package
specifically to automate some of the more tedious calculations, we are able to show that indeed
\begin{equation*}
\vartheta(n) = \frac{2}{21} + \frac{32}{441n} + \epsilon(n),
\end{equation*}
where the error term $\epsilon(n)$ is bounded in absolute value by $K n^{-7/4}$ for a constant $K$ that is computed explicitly. The precise statement and its consequences will be given later. All of the details, including links to the worksheets containing the corresponding code,
can be found in Section~\ref{sec:ramanujan-example}.

The second example we use to illustrate the methods is based on a question
from a recent paper by B\'ona and DeJonge \cite{Bona2020pattern}: let $a_n$ be the
number of 132-avoiding permutations of length~$n$ that have a unique longest 
increasing subsequence, which is also the number of plane trees with $n+1$ vertices
with a single leaf at maximum distance from the root, or the number of Dyck paths of
length $2n$ with a unique peak of maximum height. Moreover, let $p_n = \frac{a_n}{C_n}$,
where $C_n = \frac{1}{n+1} \binom{2n}{n}$ is the $n$-th Catalan number. This can
be interpreted as the probability that a 132-avoiding permutation of length $n$
chosen uniformly at random has a unique longest increasing subsequence---equivalently,
that a plane tree with $n+1$ vertices has a single leaf at maximum distance from the
root, or that a Dyck path of length $2n$ has a unique peak of maximum height.

\begin{problem}[B\'ona and DeJonge \cite{Bona2020pattern}]\label{prob:bona_dejonge}
Is it true that the sequence $p_n$ is decreasing for $n \geq 3$?
\end{problem}

While it would obviously be interesting to have a combinatorial proof, it turns out (as we will explain in the following section) that the problem can be translated in a fairly mechanical fashion (using generating functions) into a purely analytic problem: specifically, the inequality
\begin{equation}\label{eq:main_ineq}
F(n) = \sum_{k=1}^n k \sigma(k) (k^2-3n+2)(2k^2-n) \binom{2n}{n-k} < 0
\end{equation}
for all $n \geq 5$, where $\sigma(k)$ is the sum of divisors of $k$.

Again, we tackle this problem using a Mellin transform strategy. In this case,
Stirling's formula can be used in the range where $k$ is ``small'' to
find suitable approximations of the binomial coefficient $\binom{2n}{n-k}$.

As we will see, the problem is complicated in this particular instance by the occurrence of nontrivial cancellations, making precise estimates challenging. The asymptotic formula (that will be proven in this paper)
\begin{equation}\label{eq:asymptotics}
F(n) = \binom{2n}{n} \Big( - \frac{n^2}{8} + \frac{n}{24} + o(n) \Big)
\end{equation}
shows that the answer to the question of B\'ona and DeJonge is affirmative for \emph{sufficiently large} $n$. However, the $o$-notation hides an error term that is potentially huge for small values of $n$, so it is not clear what ``sufficiently large'' means.
In order to show that $p_n$ is increasing for \emph{all} $n \geq 3$, we are once again
interested in finding a quantified version of the asymptotic expansion. Again, our
SageMath package proves to be very helpful. Details are presented in
Section~\ref{sec:bona-dejonge}, including all of the code responsible
for the calculations that ultimately yield explicit error bounds for the
asymptotic expansion of $F(n)$.

\section{\texorpdfstring{$B$}{B}-terms and asymptotics with explicit error bounds}
\label{sec:bterms}

In this section, we provide the necessary background on $B$-terms and their software implementation.
We base our work on the module for computing with asymptotic expansions~\cite{Hackl-Heuberger-Krenn:2016:asy-sagemath} in SageMath~\cite{SageMath}. While
this module presently also offers some basic support for $B$-terms,
we have extended its capabilities to add support for computations
involving an additional monomially bounded ``dependent'' variable.

Before we go into the details of our new package, we briefly describe the architecture
of SageMath's asymptotic expansion module.

\subsection{Asymptotic Expansions in SageMath}

Asymptotic expansions have been introduced to SageMath in a
\emph{Google Summer of Code} project in
2015~\cite{Hackl-Heuberger-Krenn:2016:asy-sagemath} and are
part of the standard library since SageMath 6.10. As it is
custom in SageMath for mathematical objects that should support
some sort of arithmetic, they are implemented within the
\code{Element}/\code{Parent} framework. The parent structure
(which, mathematically, corresponds to the algebraic structure
containing the elements) is the so-called \code{AsymptoticRing}\footnote{
  Comprehensive documentation for the \code{AsymptoticRing} can
  be found at \url{https://doc.sagemath.org/html/en/reference/asymptotic/sage/rings/asymptotic/asymptotic_ring.html}.
}.

Broadly speaking, asymptotic expansions as implemented in SageMath
are sums of \emph{asymptotic terms} that have a certain \emph{growth}
that follows a user-specified pattern. To make this rather abstract
definition more concrete, let us consider an example.

\begin{lstlisting}[caption={Creation of an \code{AsymptoticRing}.}, label=lst:asy-creation]
 sage: A.<n> = AsymptoticRing('n^QQ * log(n)^QQ', coefficient_ring=SR, default_prec=5)
 sage: A
 Asymptotic Ring <n^QQ * log(n)^QQ> over Rational Field
 sage: A.an_element()
 1/8*n^(3/2)*log(n)^(3/2) + O(n^(1/2)*log(n)^(1/2))
\end{lstlisting}

Listing~\ref{lst:asy-creation} demonstrates how an \code{AsymptoticRing}
is created. The first argument to the constructor, the string \code{'n^QQ * log(n)^QQ'},
is a short-hand notation that specifies the ring's \emph{growth group}---the parent
that holds the ``asymptotic growth templates'' of the terms in our expansions.
The growth group specified here comprises all terms of the form $n^a \log(n)^b$,
where $a$ and $b$ are rational numbers. The second argument,
\code{coefficient_ring=SR}, specifies that the coefficients of our terms
are symbolic expressions living in SageMath's \code{SymbolicRing}.

Asking the \code{AsymptoticRing} to return an element (as in the
last line of Listing~\ref{lst:asy-creation}) yields the expansion
\[ \frac{1}{8} n^{3/2} \log(n)^{3/2} + O(n^{1/2} \log(n)^{1/2}), \]
which contains two of the three term types that are currently supported:
\emph{exact} terms (like the first summand, some non-zero constant from
the coefficient ring times a growth element) and $O$-terms. The third
type, $B$-terms, are ``quantified'' $O$-terms and will be discussed
at length in Section~\ref{sec:dependent-bterms}.

Arithmetic with asymptotic expansions works exactly as one would expect:
we can add, subtract, multiply arbitrary expansions; $O$-terms automatically
absorb terms of weaker growth. Division, raising to
a power, and taking the logarithm are all supported too---however, as the
resulting expansions are required to be a finite sum of terms, the package
computes expansions up to the specified expansion precision
(\code{default_prec=5} in Listing~\ref{lst:asy-creation}).
Examples of this mechanism are given in Listing~\ref{lst:asy-arith}.

\begin{lstlisting}[caption={Basic arithmetic with asymptotic expansions.}, label=lst:asy-arith]
 sage: n^3 + n^2 + n + O(n^2)
 n^3 + O(n^2)
 sage: 1/(n-1)
 n^(-1) + n^(-2) + n^(-3) + n^(-4) + n^(-5) + O(n^(-6))
 sage: log(1 + 1/n)
 n^(-1) - 1/2*n^(-2) + 1/3*n^(-3) - 1/4*n^(-4) + 1/5*n^(-5) + O(n^(-6))
 sage: (1 + pi/n + O(n^(-5)))^n
 e^pi - 1/2*pi^2*e^pi*n^(-1) + (1/24*(3*pi^4 + 8*pi^3)*e^pi)*n^(-2) + (-1/48*(pi^6 + 8*pi^5 + 12*pi^4)*e^pi)*n^(-3) + O(n^(-4))
\end{lstlisting}

A technical detail: as implicitly demonstrated by the last computation in
Listing~\ref{lst:asy-arith}, the \code{AsymptoticRing} supports \emph{coercion}:
as long as SageMath is able to determine a suitable inclusion of the given
parent into a compatible larger one, it will lift the elements to the larger
structure and perform the arithmetic there. In this case, the coefficient ring
of the given parent has been lifted from the rational field to the ring of symbolic
expressions in order to construct the term $\pi / n$ in our input.

Algorithmically, arithmetic with asymptotic expansions is implemented by means
of operations on
partially ordered sets, which is what expansions use as their internal data structures.
The partial order is determined by the specified growth group
(which can be significantly more complicated than a simple lexicographic order,
in particular when the growth includes exponential elements or multiple
independent variables).

\begin{lstlisting}[caption={Visualization of the partially ordered set data structure holding the terms of an asymptotic expansion.}, label=lst:asy-poset]
 sage: B.<x, y> = AsymptoticRing('x^QQ * log(x)^QQ * (QQ_+)^y * y^QQ', QQ); B
 Asymptotic Ring <x^QQ * log(x)^QQ * QQ^y * y^QQ> over Rational Field
 sage: expr = x^2 * y + x * y^2 + 42^(-y) * x + O(y / x); expr
 x*y^2 + x^2*y + O(x^(-1)*y) + x*(1/42)^y
 sage: print(expr.summands.repr_full(reverse=True))
 poset(x*y^2, x^2*y, O(x^(-1)*y), x*(1/42)^y)
 +-- oo
 |   +-- no successors
 |   +-- predecessors:   x*y^2, x^2*y
 +-- x*y^2
 |   +-- successors:   oo
 |   +-- predecessors:   O(x^(-1)*y), x*(1/42)^y
 +-- x^2*y
 |   +-- successors:   oo
 |   +-- predecessors:   O(x^(-1)*y), x*(1/42)^y
 +-- O(x^(-1)*y)
 |   +-- successors:   x*y^2, x^2*y
 |   +-- predecessors:   null
 +-- x*(1/42)^y
 |   +-- successors:   x*y^2, x^2*y
 |   +-- predecessors:   null
 +-- null
 |   +-- successors:   O(x^(-1)*y), x*(1/42)^y
 |   +-- no predecessors
\end{lstlisting}

See Listing~\ref{lst:asy-poset} for an example
where the underlying growth group is not totally ordered. In particular the
growth values of the terms $O(y/x)$ and $42^{-y} x$ are not comparable, which
is reflected in the representation of the poset where both \code{O(x^(-1)*y)}
and \code{x*(1/42)^y} have the artificial minimal element \code{null} as their
predecessor.

\subsection{The \code{dependent_bterms} package}
\label{sec:dependent-bterms}

Let us now turn to the new features that we have implemented in the \code{dependent_bterms} package. The main difference in contrast to the implementation
of $B$-terms shipped with SageMath is that we have added support for computations
with a ``dependent variable'' $k = k(n)$ that is bounded by
$c n^{\alpha} \leq k\leq C n^{\beta}$ for positive constants $c, C > 0$
and real numbers $\alpha \geq 0$ and $\beta > \alpha$.
Moreover, the package also contains several useful utility functions like,
for example, Taylor expansions with explicit error bounds.

Our package is not (yet) included in SageMath directly, but can easily be
made available to your local installation by running\footnote{
  To get the latest version of the package, in particular if an older
  version is already installed, the \code{-U} flag (for upgrade) should be
  added; i.e., \code{sage -pip install -U dependent_bterms}
  or \texttt{\%}\code{pip install -U dependent_bterms}.
}
\begin{center}
      \code{sage -pip install dependent_bterms}
\end{center}
from your terminal. Alternatively, the module can be installed by executing a cell containing
\begin{center}
      \texttt{\%}\code{pip install dependent_bterms}
\end{center}
from within a SageMath Jupyter notebook. The source code is publicly
available under a GPLv3 license on GitHub at
\url{https://github.com/behackl/dependent_bterms}. The computations in this
paper were carried out using version \code{v1.0.2} of the package.

We will now briefly walk through the core functionalities of our toolbox.
The central interface is the function
\begin{center}
\code{AsymptoticRingWithDependentVariable},
\end{center}
which generates a suitable parent structure for our desired asymptotic expansions. Listing~\ref{lst:setup} demonstrates how it is used
to instantiate the structure that will be used throughout the following
examples. We consider $1 = n^0 \leq k \leq n^{4/7}$, i.e., $\alpha = 0$
and $\beta = 4/7$.

\begin{lstlisting}[caption={Setup of the modified \code{AsymptoticRing}.}, label=lst:setup]
sage: import dependent_bterms as dbt
sage: AR, n, k = dbt.AsymptoticRingWithDependentVariable(
....:     'n^QQ', 'k', 0, 4/7, bterm_round_to=3, default_prec=5
....: )
sage: AR
Asymptotic Ring <n^QQ> over Symbolic Ring
\end{lstlisting}

\noindent The arguments passed to the interface are, in order,
\begin{itemize}
\item \code{growth_group} -- the (univariate) growth group\footnote{
      See \href{https://doc.sagemath.org/html/en/reference/asymptotic/sage/rings/asymptotic/asymptotic_ring.html}{SageMath's documentation on asymptotic expansions and the \code{AsymptoticRing}}
      for an introduction to the algebraic terminology used here.
      }
      modeling the desired structure of the asymptotic expansions.
      For example, \code{'n^QQ'} represents terms
      like $42 n^{9/13}$ or $O(n^{5/42})$.
\item \code{dependent_variable} -- a string representation of the symbolic
  variable being endowed with asymptotic growth bounds, e.g., \code{'k'}.
\item \code{lower_bound_power} -- a real number $\alpha \geq 0$ representing
  the power to which the ring's independent variable is raised when
  constructing the lower bound $c n^{\alpha} \leq k$.
\item \code{lower_bound_factor} -- a positive real number $c > 0$,
  part of the lower bound of the dependent variable. By default, $c = 1$.
\item \code{upper_bound_power} -- a real number $\beta > \alpha \geq 0$, 
  analogous to \code{lower_bound_power}, just for the upper bound
  $k \leq C n^{\beta}$.
\item \code{upper_bound_factor} -- a positive real number $C > 0$,
  part of the upper bound of the dependent variable. By default, $C = 1$.
\item \code{bterm_round_to} -- a non-negative integer or \code{None} (the default),
  determining the number of floating point digits to which the coefficients of
  $B$-terms are automatically rounded. If \code{None}, no rounding is performed.
\end{itemize}
Any other keyword arguments (like \code{default_prec} in
Listing~\ref{lst:setup} above) are passed to the constructor of
\code{AsymptoticRing}.

In this structure, arithmetic with asymptotic expansions in $n$
can be carried out in the same way as in the standard \code{AsymptoticRing},
see Listing~\ref{lst:default-arith} for a few examples that match the
behavior discussed in the previous section.

\begin{lstlisting}[caption={Arithmetic and automatic expansions in \code{AsymptoticRing}.}, label=lst:default-arith]
sage: (1 + 3*n) * (4*n^(-7/3) + 42/n + 1)
3*n + 127 + 42*n^(-1) + 12*n^(-4/3) + 4*n^(-7/3)
sage: prod((1 + n^(-j)) for j in srange(1, 10)) * (1 + O(n^(-10)))
1 + n^(-1) + n^(-2) + 2*n^(-3) + 2*n^(-4) + 3*n^(-5) + 4*n^(-6)
+ 5*n^(-7) + 6*n^(-8) + 8*n^(-9) + O(n^(-10))
sage: n / (n - 1)
1 + n^(-1) + n^(-2) + n^(-3) + n^(-4) + O(n^(-5))
\end{lstlisting}

In the implementation of the \code{AsymptoticRing} shipped with SageMath,
asymptotic expansions internally rely on ordering their summands with
respect to the growth of the independent variable(s), regardless of attached coefficients.

In the extension of our \code{dependent_bterms} module, expansions are aware
of the growth range contributed by the dependent variable appearing
in coefficients. In fact, in our modified ring, expansions are ordered with
respect to the upper bound of the coefficient growth combined with the growth
of the independent variable. This explains the---at first glance counterintuitive---ordering
of the summands in Listing~\ref{lst:dependent-arith}. The individual growth
ranges of the summands are printed at the end of the listing.

\begin{lstlisting}[caption={Arithmetic involving the dependent variable.}, label=lst:dependent-arith]
sage: k*n^2 + O(n^(3/2)) + k^3*n 
k^3*n + k*n^2 + O(n^(3/2))
sage: for summand in expr.summands.elements_topological():
....:     print(f"{summand} -> {summand.dependent_growth_range()}")
O(n^(3/2)) -> (n^(3/2), n^(3/2))
k*n^2 -> (n^2, n^(18/7))
k^3*n -> (n, n^(19/7))
\end{lstlisting}

Automatic power series expansion (with an $O$-term error) also works natively
in our modified ring, see Listing~\ref{lst:dependent-auto}. Observe that 
the error term $O(n^{-15/7})$ would actually be able to \emph{partially absorb}
some of the terms in the automatic expansion like $(k/2 + 1/6)n^{-3}$. This partial
absorption is, however, not carried out automatically due to performance reasons.
Using the \code{simplify_expansion} function included in our module expands
the symbolic coefficients and enables the error terms to carry out all allowed
(partial) absorptions.

\begin{lstlisting}[caption={Automatic expansions and manual simplifications.}, label=lst:dependent-auto]
sage: auto_expansion = exp((1 + k)/n)
sage: auto_expansion
1 + (k + 1)*n^(-1) + (1/2*(k + 1)^2)*n^(-2) + (1/6*(k + 1)^3)*n^(-3)
+ (1/24*(k + 1)^4)*n^(-4) + O(n^(-15/7))
sage: dbt.simplify_expansion(auto_expansion)
1 + (k + 1)*n^(-1) + (1/2*k^2 + k + 1/2)*n^(-2)
+ (1/6*k^3 + 1/2*k^2)*n^(-3) + 1/24*k^4*n^(-4) + O(n^(-15/7))
\end{lstlisting}

Now let us turn to the core feature of our extension: $B$-terms. In a nutshell,
$B$-terms are $O$-terms that come with an explicitly specified constant and
validity point. For example, the term $B_{n\geq 10}(42 n^3)$ represents an error term
that is bounded by $42 n^3$ for $n \geq 10$.

Listing~\ref{lst:bterms-arith} demonstrates basic arithmetic with $B$-terms. It is
worth spending a moment to understand how the resulting constants are determined.
In the first example, the $B$-term $B_{n\geq 10}(5/n)$ absorbs the exact term $3/n^2$
of weaker growth. It does so by automatically estimating $\frac{3}{n^2} \leq \frac{3}{10n}$
(as the term is valid for $n \geq 10$) and then directly absorbing the upper bound; 
$\frac{53}{10} = 5 + \frac{3}{10}$.

The same mechanism happens in the second example.
In order to avoid the (otherwise rapid) accumulation of complicated symbolic expressions
in the automatic estimates, we have specified (via the \code{bterm_round_to}-parameter that
we have set to 3) that $B$-terms should automatically be rounded to three floating point digits.
This is why the constant is given as 
$\lceil (1 + 10^{-1/3}) \cdot 10^3 \rceil \cdot 10^{-3} = \frac{293}{200}$.

\begin{lstlisting}[caption={Arithmetic with $B$-terms and the dependent variable.}, label=lst:bterms-arith]
sage: 7*n + AR.B(5/n, valid_from=10) + 3/n^2
7*n + B(53/10*n^(-1), n >= 10)
sage: AR.B(1/n, valid_from=10) + n^(-4/3)
B(293/200*n^(-1), n >= 10)
sage: AR.B(3*k^2/n^3, valid_from=10) + (1 - 2*k + 3*k^2 - 4*k^3)/n^5
B(3373/1000*abs(k^2)*n^(-3), n >= 10)
\end{lstlisting}

The third example in Listing~\ref{lst:bterms-arith} illustrates arithmetic involving
the dependent variable, which requires additional care. With $1 \leq k \leq n^{4/7}$
in place, the growth of the given $B$-term ranges from $\Theta(n^{-3})$ to
$\Theta(n^{-13/7})$. The growth of the explicit term that is added ranges from $\Theta(n^{-5})$
to $\Theta(n^{-23/7})$. In this setting, we consider the explicit term to be of
weaker growth, as the lower bound of the $B$-term is stronger than the lower bound
of the explicit term, and likewise for the upper bound. Thus we may let the $B$-term absorb it.
We do so by first estimating
\[ 
  \bigg\lvert\frac{1 - 2k + 3k - 4k^3}{n^5}\bigg\rvert
  \leq \frac{(1 + 2 + 3 + 4) k ^3}{n^5} 
  = \frac{10 k^3}{n^5}.
\]
As the power of $k$ in this bound is larger than the maximal power of $k$ in the $B$-term,
we may not yet proceed as above (otherwise we would increase the upper bound of the $B$-term,
which we must avoid). Instead, we first use $k\leq n^{4/7}$, followed by
$n\geq 10$, to obtain
\[ 
  \frac{10 k^3}{n^5} \leq \frac{10 k^2 n^{4/7}}{n^5} = \frac{10 k^2}{n^{31/7}}
  \leq \frac{10 k^2}{10^{10/7} \cdot n^3} = 10^{-3/7} \cdot \frac{k^2}{n^3},
\]
which the $B$-term can now absorb directly. Hence the value of the $B$-term constant
is determined by $\lceil (3 + 10^{-3/7}) \cdot 10^3 \rceil \cdot 10^{-3} = \frac{3373}{1000}$.

Finally, our module also provides support for $B$-term bounded Taylor expansion (again,
also involving the dependent variable) in form of the \code{taylor_with_explicit_error}
function. An example is given in Listing~\ref{lst:bterm-taylor}: we first obtain a Taylor
expansion of $f(t) = (1 - t^2)^{-1}$ around $t = (1 + k)/n + B_{n\geq 10}(k^3/n^3)$.
Using the \code{simplify_expansion} function rearranges the terms and
lets the $B$-term absorb coefficients (partially) as far it is able to.
Observe that it may happen
that the attempted simplification produces summands with a smaller upper growth
bound that the implementation cannot absorb ($B_{n\geq 10}(k^3/n^3)$ vs. $n^{-2}$
in this case). The expansion is still correct; just not as compact as it could be.
We can also use the \code{simplify_expansion} function with the \code{simplify_bterm_growth}
parameter set to \code{True} to collapse the dependent
variables in all $B$-terms by replacing them with their upper bounds, resulting in
a single ``absolute'' $B$-term.

\begin{lstlisting}[caption={$B$-term bounded Taylor expansions.}, label=lst:bterm-taylor]
sage: arg = (1 + k)/n + AR.B(k^3/n^3, valid_from=10)
sage: ex = dbt.taylor_with_explicit_error(
....:     lambda t: 1/(1 - t^2), arg, 
....:     order=3, valid_from=10)
sage: ex
1 + ((k + 1)^2)*n^(-2)
+ B((abs(7351/250*k^3 + 30*k^2 + 30*k + 10))*n^(-3), n >= 10)
sage: dbt.simplify_expansion(ex)
1 + k^2*n^(-2)
+ B((abs(7351/250*k^3 + 30*k^2 + 30*k + 10))*n^(-3), n >= 10)
+ (2*k + 1)*n^(-2)
sage: dbt.simplify_expansion(ex, simplify_bterm_growth=True)
1 + k^2*n^(-2) + B(41441/1000*n^(-9/7), n >= 10)
\end{lstlisting}

\section{A detailed example: the Ramanujan \texorpdfstring{$Q$}{Q}-function}
\label{sec:ramanujan-example}

As announced in the introduction, we now demonstrate how the
\code{dependent_bterms} package can be used in practice. Recall
that Problem~\ref{prob:ramanujan}, based on Ramanujan's Question 294,
asks to find a precise asymptotic expansion for the difference between
two functions, $Q(n)$ and $R(n)$, that are given by
\[
Q(n) = \sum_{k=0}^{n-1} \frac{n!}{(n-k-1)!n^{k+1}}
= 1 + \frac{n-1}{n} + \frac{(n-1)(n-2)}{n^2} + \cdots + \frac{(n-1)!}{n^{n-1}}
\]
and
\[
R(n) = \sum_{k=0}^{\infty} \frac{n!n^k}{(n+k)!}
= 1 + \frac{n}{n+1} + \frac{n^2}{(n+1)(n+2)} + \cdots.
\]
To be more precise, we want to show that
\begin{equation}\label{eq:ramanujan-theta}
R(n) - Q(n) = \frac23 + \frac{8}{135(n + \vartheta)}
\end{equation}
for some $\vartheta = \vartheta(n)$ between $\frac{2}{21}$ and $\frac{8}{45}$.
From works of Flajolet, Grabner, Kirschenhofer, and Prodinger~\cite{Flajolet-Grabner-Kirschenhofer-Prodinger:1995:ramanujan-q},
Ramanujan himself, as well as Watson~\cite{Watson:1929:theorems-ramanujan} and Knuth~\cite{Knuth:1973:TAOCPv3},
we know that $Q(n)$ and $R(n)$ admit the asymptotic expansions
\[
Q(n) \sim \sqrt{\frac{\pi n}{2}} - \frac13
+ \frac1{12} \sqrt{\frac{\pi}{2n}} - \frac{4}{135n} + \cdots
\]
and
\[
R(n) \sim \sqrt{\frac{\pi n}{2}} + \frac13 + \frac1{12} \sqrt{\frac{\pi}{2n}}
+ \frac{4}{135n} + \cdots.
\]
The proofs are based on integral representations: specifically, $Q(n)$ can be expressed as
\[ Q(n) = \int_0^{\infty} e^{-x} \Big(1 + \frac{x}{n} \Big)^{n-1}\,dx. \]
Flajolet, Grabner, Kirschenhofer and
Prodinger~\cite{Flajolet-Grabner-Kirschenhofer-Prodinger:1995:ramanujan-q},
on the other hand, use a complex integral, namely
\[
R(n) - Q(n) =
\frac{n!}{n^{n-1}} \frac{1}{2\pi i} \oint_{\mathcal{C}}
\log \frac{(1-y)^2}{2(1-y e^{1-y})} (1-y)e^{ny} \frac{dy}{y^{n+1}}
\]
for a small contour $\mathcal{C}$ around the origin.

Here we apply a ``brute force'' approach based directly on the series representations
that does not require any of these ingenious integral representations. To this end,
we split the sum
\[
R(n) - Q(n) =
\sum_{k=1}^{\infty} \Big(\frac{n!n^k}{(n+k)!} - \frac{n!}{(n-k-1)!n^{k+1}} \Big),
\]
where we interpret the second fraction as $0$ for $k \geq n$, into parts. To be more
precise, we choose parameters $\alpha, C > 0$ and then consider three cases for $k$:
The \emph{small} range, where $1 \leq k < C n^{\alpha}$, the \emph{middle} range
where $C n^{\alpha} < k < n$, and the \emph{large} range where $k \geq n$.

Our strategy for treating the asymptotic contribution of these sums over the
individual ranges now consists of first checking that the quantity $\vartheta(n)$
from~\eqref{eq:ramanujan-theta}
indeed satisfies $\frac{2}{21} \leq \vartheta(n) \leq \frac{8}{45}$ for all $n$
up to some threshold value, $n \leq N = 70000$. This allows us to consider
$n\geq N$ when constructing $B$-terms that hold the quantified expansion
errors. All computations can be found in the ancillary SageMath notebook
located at
\begin{center}
  \url{https://arxiv.org/src/2403.09408v3/anc/2025-ramanujan-q.ipynb},
\end{center}
and for the sake of convenience a static version of this notebook
showing the executed results of all computations can be found at
\begin{center}
  \url{https://arxiv.org/src/2403.09408v3/anc/2025-ramanujan-q.html}.
\end{center}
We will, however, for the sake of demonstration, also include individual
code snippets that illustrate how the \code{dependent_bterms} package is used
throughout this section.

Now, before we start estimating the individual contributions, let us first
observe some additional useful facts. First, note that when studying the
difference $R(n) - Q(n)$, all individual differences in the sum
\[
  R(n) - Q(n) = \sum_{k\geq 1} \Bigl( \frac{n!\, n^k}{(n + k)!} - \frac{n!}{(n-k-1)!\, n^{k+1}} \Bigr)
\]
are non-negative as we can write
\begin{align*}
  \frac{n!\, n^k}{(n+k)!}
  &= \frac{n^k}{(n+1) (n+2) \cdots (n+k)}
  = \prod_{j=1}^k \Bigl( 1 + \frac{j}{n} \Bigr)^{-1} \\
  &\geq \prod_{j=1}^k \Bigl( 1 - \frac{j}{n} \Bigr)
  = \frac{(n-1)(n-2) \cdots (n-k)}{n^k} = \frac{n!}{(n-k-1)!\, n^{k+1}}.
\end{align*}
This allows us to focus on deriving upper bounds only.

Secondly, we are able to determine a nice exponential bound for the
summands in $Q(n)$ by rewriting the product representation above as a
sum of logarithms, and then using that $\log(1+x) \geq x - x^2 / 2$ for $x \geq 0$ (it is easy to verify that the derivative of $\log(1+x) - x + x^2 / 2$ is $\frac{x^2}{1+x} \geq 0$ while its value at $0$ is $0$). We find
\begin{multline}
  \frac{n! n^k}{(n+k)!} = \exp\Bigl( -\sum_{j=1}^k \log(1 + j/n) \Bigr) \\
  \leq \exp\Bigl( -\frac{1}{n} \sum_{j=1}^k j + \frac{1}{2n^2} \sum_{j=1}^k j^2 \Bigr)
  = \exp\Bigl( -\frac{(6n - 2k - 1) (k+1) k}{12n^2} \Bigr),
\end{multline}\label{eq:ramanujan:exponential}
which will be useful in two out of our three respective summation ranges.

\subsection{Contribution of the large range}
Let us consider the contribution of the tail of the summation, the range
where $k \geq n$. Observe that, as mentioned above, the summands in $Q(n)$
collapse to 0 in this range. Furthermore, we have the inequality
$(n+k)! \geq (2n)!\, (2n)^{k-n}$, which then yields
\begin{align}
  \sum_{k\geq n} \frac{n!\, n^k}{(n + k)!}
  &\leq \sum_{k\geq n} \frac{n!\, n^k}{(2n)! (2n)^{k-n}}
  = \frac{n!\, n^n}{(2n)!} \sum_{k \geq n} \frac{1}{2^{k-n}}
  = 2 \frac{n!\, n^n}{(2n)!} \notag \\
  & \leq 2 \exp\Bigl( -\frac{(4n-1)(n+1)n}{12 n^2} \Bigr) \leq 2e^{-n/3},
  \label{eq:ramanujan:error-large}
\end{align}
an exponentially small error.

\subsection{Contribution of the middle range}
Consider $C n^{\alpha} < k < n$. If we choose the ``cutoff'' $C n^{\alpha}$
sufficiently large, we can obtain a suitable error term by estimating the individual
differences by first ignoring the contribution of $Q(n)$, and then working with
the exponential bound~\eqref{eq:ramanujan:exponential}. That is,
\begin{align*}
  \frac{n! n^k}{(n+k)!} - \frac{n!}{(n-k-1)! n^{k+1}}
  &\leq \frac{n! n^k}{(n+k)!}
  \leq \exp\Bigl( -\frac{(6n - 2k - 1) (k+1) k}{12n^2} \Bigr) \\
  &\leq \exp\Bigl( -\frac{k^2}{2n} + \frac{k^3}{6n^2} \Bigr).
\end{align*}
This estimate is increasing for all $1\leq k \leq 2n$, which means that we
can bound the sum over the range $Cn^{\alpha} < k < n$ by means of an integral,
\begin{align*}
  \sum_{Cn^{\alpha} < k < n} \exp\Bigl( -\frac{k^2}{2n} + \frac{k^3}{6n^2} \Bigr)
  &\leq \int_{C n^{\alpha}}^n \exp\Bigl(-\frac{t^2}{2n} + \frac{t^3}{6n^2}\Bigr) \mathrm{d}t
  \\
  & \leq \int_{C n^{\alpha}}^n 
  \frac{\frac{t}{n} - \frac{t^2}{2n^2}}{C n^{\alpha - 1} (1 - C n^{\alpha - 1} / 2)}
  \exp\Bigl(-\frac{t^2}{2n} + \frac{t^3}{6n^2}\Bigr) \mathrm{d}t \\
  &= 
  \frac{1}{C n^{\alpha - 1} (1 - C n^{\alpha - 1} / 2)}
  \Bigl[ - \exp\Bigl( -\frac{t^2}{2n} + \frac{t^3}{6n^2} \Bigr) \Bigr]_{C n^{\alpha}}^n \\
  &= 
  \frac{
    \exp(- C^2 n^{2\alpha - 1} / 2 + C^3 n^{3\alpha - 2} / 6) - \exp(-n/3)
  }{C n^{\alpha - 1} (1 - C n^{\alpha - 1} / 2)}.
\end{align*}

In order to ensure that the contribution of the sum over this range is
exponentially small, our parameter $\alpha$ needs to be chosen such that
$2\alpha - 1$ is positive, while $3\alpha - 2$ is negative. In other words,
$\alpha$ has to be chosen from the interval $(\frac12, \frac23)$.

Combining the (exponentially small) contribution of the large range with
the contribution determined above, and with some further simplifications
we arrive at
\begin{align*}
  \sum_{k > C n^{\alpha}} \Bigl( \frac{n! n^k}{(n+k)!} - \frac{n!}{(n-k-1)! n^{k+1}} \Bigr)
  &\leq
  \frac{n^{1 - \alpha}}{C} \frac{1}{1 - \frac{C}{2} n^{\alpha - 1}}
  \exp\Bigl( -\frac{C^2}{2} n^{2\alpha - 1} + \frac{C^3}{6} n^{3\alpha - 2} \Bigr) \\
  &\leq
  \frac{2 \exp(C^3 N^{3\alpha - 2} / 6)}{C} n^{1 - \alpha} e^{- \frac{C^2}{2} n^{2\alpha - 1}}.
\end{align*}
The role of the parameter $C > 0$ is now to accelerate the exponential
decay of the contribution over this range such that it is numerically negligible
for $n > N$. For our purposes, choosing $\alpha = \frac{3}{5}$ and $C = 4$ is
appropriate. A crude estimate then yields the comparatively very small
$B$-term
\begin{equation}\label{eq:ramanujan:tails-pruning}
  \frac{2 \exp(C^3 N^{3\alpha - 2} / 6)}{C} n^{1 - \alpha} e^{- \frac{C^2}{2} n^{2\alpha - 1}}
  = B_{n\geq N}(3.11 \cdot 10^{-31} n^{-10})
\end{equation}
as a bound of the contribution over the medium and large range of the sum.
Using our package, this bound can be derived by following the steps outlined in
Listing~\ref{lst:ramanujan-tail-error}.

\begin{lstlisting}[caption={
  Finding an explicit bound for the tail error. The variables
  \code{n_sym} and \code{n} represent $n$ in the \code{SymbolicRing}
  and a suitable \code{AsymptoticRing}, respectively.
}, label=lst:ramanujan-tail-error]
 sage: tail_error_sym = 2 / C * exp(C^3 * N^(3*alpha - 2) / 6) \
 ....:   * n_sym^(1 - alpha) * exp(-C^2 / 2 * n_sym^(2*alpha - 1))
 sage: bool((n_sym^10 * tail_error_sym).diff(n_sym, 1)(n_sym=N) < 0)  # ensure N is large enough such that bound is decreasing
 True
 sage: tail_error_asy = tail_error_sym(m=N) \
 ....:   * AR.B(n^(-10), valid_from=N)
 sage: tail_error_asy
 B(3.105655...e-31*n^(-10), n >= 70000)
\end{lstlisting}

\subsection{Contribution of the small range}
The last case, the range $1 \leq k < C n^{\alpha}$, carries the main contribution
of the sum. To analyze this part, we first find a suitable expansion
of the summands that allows rewriting the sum in such a way that allows
using the Mellin transform to determine the asymptotic growth.

Observe that the product representations
of the two terms in the difference can be used to obtain
\begin{align*}
  \frac{n! n^k}{(n+k)!}
  &= \prod_{j=1}^k \Bigl(1 + \frac{j}{n}\Bigr)^{-1}
  = \exp\Bigl(- \sum_{j=1}^k \log(1 + j/n)\Bigr)
  = \exp\Bigl(\sum_{j=1}^k \sum_{r \geq 1} \frac{(-1)^r j^r}{r n^r}\Bigr),
  \\
  \frac{n!}{(n-k-1)! n^{k+1}}
  &= \prod_{j=1}^k \Bigl(1 - \frac{j}{n}\Bigr)
  = \exp\Bigl( \sum_{j=1}^k \log(1 - j/n) \Bigr)
  = \exp\Bigl( -\sum_{j=1}^k \sum_{r \geq 1} \frac{j^r}{r n^r} \Bigr).
\end{align*}
Hence, the difference can be expressed as
\begin{align*}
\frac{n! n^k}{(n+k)!} - \frac{n!}{(n-k-1)! n^{k+1}}
&=
\exp\Bigl(\sum_{j=1}^k \sum_{r \geq 1} \frac{(-1)^r j^r}{r n^r}\Bigr)
-
\exp\Bigl( -\sum_{j=1}^k \sum_{r \geq 1} \frac{j^r}{r n^r} \Bigr) \\
&= 
\exp\Bigl(- \sum_{j=1}^k \sum_{\substack{r \geq 1 \\ r \text{ odd}}}
  \frac{j^r}{r n^r}\Bigr) \\
  & \qquad \times \biggl(
  \exp\Bigl(\sum_{j=1}^k \sum_{\substack{r \geq 1 \\ r \text{ even}}}
  \frac{j^r}{r n^r}\Bigr) - 
  \exp\Bigl(- \sum_{j=1}^k \sum_{\substack{r \geq 1 \\ r \text{ even}}}
  \frac{j^r}{r n^r}\Bigr)
\biggr)\\
&= 2 \exp\Bigl(- \sum_{j=1}^k \sum_{\substack{r \geq 1 \\ r \text{ odd}}}
  \frac{j^r}{r n^r}\Bigr)
  \sinh\Bigl(\sum_{j=1}^k \sum_{\substack{r \geq 1 \\ r \text{ even}}}
  \frac{j^r}{r n^r}\Bigr).
\end{align*}

To find a suitable (bounded) asymptotic expansion from this representation,
we choose to cut the sums appearing as arguments of $\exp$ and $\sinh$
at some point (practically, we choose $R = 9$) and then estimate them
by means of

\begin{equation*}
0 \leq \sum_{j=1}^k \sum_{\substack{ r=R \\ r \text{ odd}}}^{\infty} \frac{j^r}{rn^r} \leq \sum_{j=1}^k \frac{j^R}{Rn^R(1-j^2/n^2)}
\end{equation*}
by means of a geometric series (observe that the factor $1-j^2/n^2$ in the denominator stems from the fact that we are only summing over odd $r$). An analogous estimate also holds for the sum
over even $r$. Furthermore, we find
\begin{align} 
\sum_{j=1}^k \frac{j^R}{Rn^R(1-j^2/n^2)}
&\leq \frac{1}{Rn^R(1-k^2/n^2)} \Big( k^R + \int_0^k t^R \,dt \Big) \nonumber \\
&= \frac{1}{Rn^R(1-k^2/n^2)} \Big( k^R + \frac{k^{R+1}}{R+1} \Big),\label{eq:sumestimate}
\end{align}
which we can use to construct a corresponding $B$-term using our software package.

\begin{lstlisting}[caption={
  Computing an explicit expansion of $2\exp(\cdots)$.
  }, label=lst:ramanujan:odd-sum-exp]
 sage: sum_odd_error = (dbt.taylor_with_explicit_error(lambda t: 1/(1 - t), k^2/n^2, valid_from=N)
 ....:   * AR.B(1 / (R * n^R) * (k^R + k^(R+1)/(R+1)), valid_from=N))
 sage: sum_odd_expansion = -sum([sum(j^r, j, 1, k) / (r * n^r) for r in srange(1, R, 2)])
 sage: sum_odd_expansion += k^2/(2*n)  # remove asymptotic main term
 sage: sum_odd_expansion += sum_odd_error
 sage: exp_expansion = dbt.taylor_with_explicit_error(lambda t: 2*exp(t), sum_odd_expansion, order=taylor_order, valid_from=N)
 sage: exp_expansion
 2 - k*n^(-1) + (-1/6*k^4 - 1/3*k^3 - 1/6*k^2)*n^(-3) + ... + B(...*k^10*n^(-9), n >= 70000) + ...
\end{lstlisting}

The computation of the corresponding Taylor expansion with quantified error
is carried out automatically using the \code{dependent_bterms} package; details are given
in Listing~\ref{lst:ramanujan:odd-sum-exp}. Combining this expansion with the expansion
obtained from analogous computations for the expansion of $\sinh(\cdots)$
\begin{align*}
  2 &\exp \Big(
    {-}\sum_{j=1}^k \sum_{\substack{r=1 \\ r \text{ odd}}}^{\infty} \frac{j^r}{r n^r}
  \Big) \sinh \Big(
    \sum_{j=1}^k \sum_{\substack{r=1 \\ r \text{ even}}}^{\infty} \frac{j^r}{r n^r}
  \Big) = 
  e^{-k^2/(2n)} \Big(
    \frac{2k^3 + 3k^2 + k}{6n^2} - \frac{2k^4 + 3k^3 + k^2}{12n^3} \\
    & + \frac{34k^5 + 75k^4 + 45k^3 - 4k}{240n^4}
    - \frac{40k^7 + 222k^6 + 375k^5 + 225k^4 + 20k^3 - 12k^2}{1440 n^5} + \cdots \\
    & + B_{n\geq N}\Bigl( \frac{730 k^{33} + 91 k^{32}}{1000000000}  \Bigr) + \cdots
  \Big),
\end{align*}  
where $\cdots$ denotes omitted terms (as the total expression, a expansion of terms
going up to $k^{58} / n^{42}$, is very long). The full expansion can be found
in the corresponding SageMath notebook at \url{https://arxiv.org/src/2403.09408v3/anc/2025-ramanujan-q.ipynb}.

Qualitatively, if we collect all ``exact'' terms in this sum in a polynomial in $k$
and $n^{-1}$ and denote this part as $S(n, k)$, and if we denote the sum of all
$B$-terms as $S_B(n, k)$, then our sum of interest can be expressed as
\begin{align*}
  \sum_{1 \leq k < C n^{\alpha}} e^{-k^2 / (2n)} (S(n, k) + S_B(n, k)),
\end{align*}
which we now want to study using the Mellin transform (see~\cite{Flajolet-Gourdon-Dumas:1995:mellin} for a general reference).

Let us study the contribution of the error part. For a positive function $f(t)$ that is increasing up to some maximum $t_0$ and decreasing thereafter,
it is well known that 
\begin{equation}\label{eq:integral-inequality}
\sum_{k = 1}^K f(k) \leq f(t_0) + \int_0^{K} f(t)\,dt.
\end{equation}
For a summand
$B_{n \geq N}(c_{a, b} k^a n^b)$ in $S_B(n, k)$ we make use of the fact that
the function $f(t) = t^a e^{-t^2 / (2n)}$ is increasing from $0$ to $t_0 = \sqrt{an}$
and decreasing afterwards, allowing us to use the integral inequality
\[ \sum_{1\leq k < C n^{\alpha}} f(k) \leq f(t_0) + \int_{0}^{C n^{\alpha}} f(t)~dt. \]
Hence we can estimate
\begin{align*}
\sum_{1 \leq k < C n^{\alpha}} B_{n\geq N}(c_{a, b} k^a n^b) e^{-k^2 / (2n)}
& \leq c_{a, b} n^b \Bigl(
  (an)^{a/2} e^{-a/2} + \int_0^{C n^{\alpha}} t^a e^{-t^2 / (2n)}~dt
  \Bigr) \\
& \leq c_{a, b} n^b \Bigl(
  (an)^{a/2} e^{-a/2} + \int_0^{\infty} t^a e^{-t^2 / (2n)}~dt \Bigr) \\
& = c_{a, b} n^{a/2 + b} \Bigl((a/e)^{a/2} + \sqrt{2^{a-1} n}\, \Gamma((a+1)/2)\Bigr).
\end{align*}

Plugging in the summands from $S_B(n, k)$ and combining the individual contributions
yields an overall error of
\begin{equation}\label{eq:ramanujan:expansion-error}
  \sum_{1 \leq k < C n^{\alpha}} e^{-k^2 / (2n)} S_B(n, k)
  = B_{n\geq N}\bigl(5913825.038\, n^{-4}\bigr).
\end{equation} 

Next we want to complete the tails: we need to make sure that extending the
summation range of $S(n, k) e^{-k^2 / (2n)}$ over $k > C n^{\alpha}$ also
only adds an error that is negligible. Practically, we can do this by
considering $S(n, k) / k^{\deg_k S(n, k)}$, an expression in which only
non-positive powers of $k$ appear, then we replace $k$ with the lower bound
and let our package compute the corresponding $B$-term. Provided that $k > C n^{\alpha}$,
we find
\[
  \frac{S(n, k)}{k^{58}} = B_{n\geq N}(3.18\cdot 10^{-34}\, n^{-35}),
\]  
or equivalently
\[
  S(n, k) = B_{n\geq N}\Bigl(3.18\cdot 10^{-34}\, \frac{k^{58}}{n^{35}}\Bigr),
\]
which we can then use with the same integral bound as above to get
\begin{align}
  \notag \sum_{k > C n^{\alpha}} S(n, k) e^{-k^2 / (2n)}
  & \leq \frac{3.18\cdot 10^{-34}}{n^{35}} \sum_{k > C n^{\alpha}} k^{58} e^{-k^2 / (2n)} \\ 
  \notag & \leq \frac{3.18\cdot 10^{-34}}{n^{35}} \int_{0}^{\infty} t^{58} e^{-t^2 / (2n)}~dt \\
  \label{eq:ramanujan:tails-completion}
  & = B_{n\geq N}(197293.89 n^{-11/2}).
\end{align}  
The computations behind this estimate are carried out in
Listing~\ref{lst:ramanujan:error-tails}.

\begin{lstlisting}[caption={
  Computation of $B$-term representing the error made when completing the tails
  of the sum over $S(n, k)$.
 }, label=lst:ramanujan:error-tails]
 sage: expansion_deg_k = expansion_symbolic.degree(k); expansion_deg_k
 56
 sage: S_nk_normalized_bound = dbt.evaluate(
 ....:     (expansion_symbolic.subs(n=n_sym) / k^expansion_deg_k)(k=C*n_sym^alpha).expand(),
 ....:     n=n).B(valid_from=N)
 sage: S_nk_normalized_bound
 B(3.17900...e-34*n^(-35), n >= 70000)
 sage: f(t) = t^j * exp(-t^2 / 2 / n_sym)
 sage: S_nk_tail_error = dbt.evaluate(
 ....:     integrate(f(t, j=expansion_deg_k), t, 0, oo).expand(), n=n
 ....: ) * S_nk_normalized_bound
 sage: S_nk_tail_error
 B(197293.8823...*n^(-11/2), n >= 70000)
\end{lstlisting}

As long as we keep track of this error, we may now study the infinite sum
\[ \sum_{k\geq 1} S(n, k) e^{-k^2 / (2n)}, \]
which, given the nature of $S(n, k)$, is a linear combination of sums of the
form
\[ \sum_{k\geq 1} k^a n^b e^{-k^2 / (2n)} =: g_{a, b}(1/n). \]
A straightforward computation shows that the Mellin transform of such
a summand $g_{a, b}(t)$ is then given by
\begin{equation}
  g_{a, b}^*(s) = \int_0^{\infty} g_{a, b}(t) t^{s-1}~dt
  = 2^{s-b} \Gamma(s-b) \zeta(2s - a - 2b),
\end{equation}
with simple poles present at $s \in b + \mathbb{Z}_{< 0}$
(originating from $\Gamma$) and $s = b + \frac{a+1}{2}$
(from $\zeta$). Note that as we know $a \geq 1$ and $b \leq -1$,
the singularities contributed by $\Gamma$ and $\zeta$ never combine
to a higher-order pole.

With the Mellin transform at hand, using the converse mapping theorem
allows us to express
\[
  g_{a, b}(t) =
  \frac{1}{2\pi i} \int_{\sigma - i\infty}^{\sigma + i\infty} g_{a, b}^*(s) t^{-s}~ds
\]
for suitable $\sigma > b + \frac{a + 1}{2}$. The vertical line of integration
may be shifted to the left, as long as we collect the residues of all poles that
we cross.

While doing so by hand is a somewhat daunting task (with the chosen expansion
precision, $S(n, k)$ consists of over 900 summands), the computer can handle
the transformation effortlessly. Collecting (and checking for) residues in
all half-integers greater than or equal to $-7/2$ reveals that the sum grows
asymptotically like
\begin{equation}\label{eq:ramanujan:mellin-exact}
  \sum_{k\geq 1} S(n, k) e^{-k^2 / (2n)} = \frac{2}{3} + \frac{8}{135 n} - \frac{16}{2835 n^2} - \frac{32}{8505 n^3} + O(n^{-15/4}),
\end{equation}
with the error term still to be quantified (and originating from the shifted integrals).
Listing~\ref{lst:ramanujan:residues} shows how a symbolic representation of the
summand expansion is processed to first determine the Mellin transform, and then
to compute the sum of all residues of interest.

\begin{lstlisting}[caption={
  Computation of Mellin transform and sum of residues.
}, label=lst:ramanujan:residues]
 sage: def mellin_from_exponents(a, b):
 ....:     return 2^(s - b) * gamma(s - b) * zeta(2*s - a - 2*b)
 sage: mellin_transform = 0
 sage: for k_poly, b in expansion_symbolic.coefficients(n_sym):
 ....:     for coef, a in k_poly.coefficients(k):
 ....:         mellin_transform += coef * mellin_from_exponents(a, b)
 sage: mellin_transform *= t^(-s)
 sage: growth_exact = dbt.evaluate(
 ....:     sum(mellin_transform.residue(s==s0) for s0 in srange(-7/2, 30, 1/2)),
 ....:     t=n^(-1))
 sage: growth_exact
 2/3 + 8/135*n^(-1) - 16/2835*n^(-2) - 32/8505*n^(-3)
\end{lstlisting}

In order to find an estimate for the contribution of the shifted integral from
the converse mapping theorem we investigate, individually for each summand in the Mellin transform, how far we can shift the line of integration
to the left (in half-integer units) until $\Re(s) = c_{a,b}$ without
collecting any further residues.

In a central region of $c_{a,b} + iw$ for $|w| \leq 100$ we use rigorous
integration via interval arithmetic to determine the numerical value of the shifted
integrals. Outside, for $|w| > 100$, we determine a suitable upper bound
of the integrand. For $\Gamma(c_{a, b} + iw)$ where $c_{a,b} > 0$ 
we use \cite[\href{https://dlmf.nist.gov/5.6.E9}{(5.6.9)}]{NIST:DLMF:v1.1.12},
and when $c_{a,b} < 0$ we first shift the argument to the right
via the functional equation $\Gamma(s) = \frac{1}{s} \Gamma(s+1)$ and then proceed as before.
For $\zeta(c_{a, b} + iw)$ we bound the modulus from above by $\zeta(3/2)$
if $c_{a, b} \geq 3/2$. When $c_{a, b} \leq -1/2$ we first apply the
reflection formula~\cite[\href{https://dlmf.nist.gov/25.4.E1}{(25.4.1)}]{NIST:DLMF:v1.1.12}; the resulting factors can all be estimated directly.
For the special case of $c_{a, b} = 1/2$ we use the bound proved
by Hiary, Patel, and Yang in
\cite[Theorem 1.1]{Hiary-Patel-Yang:2024:estimate-zeta} to obtain
\[ |\zeta(1/2 + iw)| \leq 0.618\, t^{1/6}\log t \leq 0.618\, t^{1/2} \]
for $t\geq 100$.

Letting a computer collect and combine all of these estimates lets us find the $B$-term
\begin{equation}\label{eq:ramanujan:shifted-integrals}
\abs*{\sum_{a, b} \frac{1}{2\pi i} \int_{-15/4-i\infty}^{-15/4 + i\infty} g_{a,b}^*(s) t^{-s}~ds}
= B_{n\geq N}\bigl(67990.898\, n^{-15/4}\bigr)
\end{equation}  
for the contribution of the shifted integrals.

We are finally at the point where we can combine all the individual asymptotic
contributions and errors we have collected throughout this section. Summing our
findings from
\eqref{eq:ramanujan:tails-pruning},
\eqref{eq:ramanujan:expansion-error},
\eqref{eq:ramanujan:tails-completion},
\eqref{eq:ramanujan:mellin-exact},
\eqref{eq:ramanujan:shifted-integrals}
yields the asymptotic expansion
\begin{equation}\label{eq:expl-asy}
  R(n) - Q(n) = \frac{2}{3} + \frac{8}{135 n} - \frac{16}{2835 n^2} - \frac{32}{8505 n^3}
  + B_{n\geq N}(431565.85 n^{-15/4}).
\end{equation} 
Recall now Ramanujan's original statement, namely that
  \begin{equation*}
    R(n) - Q(n) = \frac{2}{3} + \frac{8}{135(n + \vartheta)}
  \end{equation*}
for some $\vartheta = \vartheta(n)$ that satisfies $2/21 \leq \vartheta(n) \leq 8/45$.
Solving for $\vartheta$, our formula~\eqref{eq:expl-asy} becomes
\begin{equation}\label{eq:ramanujan:theta-asy}
\vartheta(n) = \frac{2}{21} + \frac{32}{441n} + B_{n\geq N}(7282696.17 n^{-7/4}).
\end{equation}
This implies the following:
\begin{itemize}
\item The constant $2/21$ in Ramanujan's statement is indeed best possible and attained in the limit.
\item On the other hand, the constant $8/45$ is attained for $n=0$, and the largest
  value of $\vartheta$ for positive $n$ is in fact
  $\vartheta(1) = \frac{368-135e}{135e-360} \approx 0.148098$
  (compared to $8/45 = 0.17777\ldots$), which is proved
  by~\eqref{eq:ramanujan:theta-asy} combined with evaluating $\vartheta$ for $n \leq N$.
\end{itemize}

Note that our asymptotic expansion effectively gives a decreasing
upper bound (seen as a function in $n$) for $\vartheta$---and
while it is not quite enough to prove it, it also makes it plausible
that $\vartheta$ is decreasing itself (which we have explicitly verified for
``small'' $n$).

\section{On a Question of B\'ona and DeJonge}
\label{sec:bona-dejonge}

Let us briefly recall the setting of Problem~\ref{prob:bona_dejonge}
as outlined in the introduction. The numbers $a_n$ count all $132$-avoiding
permutations of length $n$ with a unique longest increasing subsequence. By
a standard bijection, $a_n$ also describes the number of plane trees with
$n+1$ vertices with a single leaf at maximum distance from the root---or,
equivalently, the number of Dyck paths of length $2n$
(i.e., lattice paths starting at $(0,0)$ and ending at $(2n,0)$ whose steps
are either ``up'' $(1,1)$ or ``down'' $(1,-1)$) with a unique peak of maximum height.

Such a path can be decomposed into two pieces: before and after the peak. The part before the peak needs to be a path that finishes at its maximum height $h$ (but does not reach it earlier, since the peak is unique), and the path after the peak needs to be a path that starts at its maximum height $h$ and never returns to it (which can also be seen as the reflection of a path that finishes at its maximum height but does not reach it earlier). Such paths were analyzed in \cite{Bousquet2008culminating} and \cite{Hackl2016analysis}. Specifically, \cite[Proposition 2.1]{Hackl2016analysis} states that the probability that a simple symmetric random walk of length $n$ never drops below $0$ and finishes at its maximum height $h$ (which can also be reached earlier) is precisely
\begin{equation*}
2 [z^{n+1}] \frac{1}{U_{h+1}(1/z)},
\end{equation*}
where $U_{h+1}$ is the Chebyshev polynomial of the second kind of degree $h+1$. A path that finishes at its maximum height $h$ without reaching that height before is obtained from a path that finishes at its maximum height $h-1$ by adding one more step up. Since every path of length $n$ has probability $2^{-n}$ to occur under the model of a simple symmetric random walk, it follows that the (ordinary) generating function for paths of maximum height $h$ that finish at the maximum and do not reach it earlier is 
\begin{equation*}
\sum_{n \geq 0} 2^n x^{n+1} \cdot 2 [z^{n+1}] \frac{1}{U_h(1/z)} = \frac{1}{U_h(1/(2x))}.
\end{equation*}
For example,
\begin{equation*}
\frac{1}{U_3(1/(2x))}= \frac{x^3}{1-2x^2} = x^3 + 2x^5 + 4x^7 + \cdots
\end{equation*}
is the generating function for paths that finish at their maximum height $3$ and do not reach this height before the final step. The formula is even true for $h=0$: $\frac{1}{U_0(1/(2x))} = 1$ is indeed the correct generating function in this case.

Since the paths we are interested in can be seen as pairs of paths that finish at their maximum height and do not reach this height before, we find that the generating function of $a_n$ is
\begin{equation*}
A(x) = \sum_{n=0}^{\infty} a_n x^{2n} = \sum_{h=0}^{\infty} \frac{1}{U_h(1/(2x))^2}.
\end{equation*}
We can simplify the expression by means of the substitution $x = \frac{\sqrt{t}}{1+t}$. Note that this yields $\frac{1}{2x} = \frac{1+t}{2\sqrt{t}} = \cosh(\frac12 \log t)$. Since $U_h(\cosh w) = \frac{\sinh((h+1)w)}{\sinh w}$, this implies that
\begin{equation*}
U_h(1/(2x)) = \frac{\sinh(\frac{h+1}{2} \log t)}{\sinh (\frac12 \log t)} = t^{-h/2} \cdot \frac{1-t^{h+1}}{1-t}.
\end{equation*}
Thus
\begin{equation*}
A(x) = \sum_{h=0}^{\infty} \frac{t^h(1-t)^2}{(1-t^{h+1})^2}.
\end{equation*}
Now we can obtain an alternative expression for $a_n$ by applying Cauchy's integral formula to the generating function $A(x)$. For suitable contours $\mathcal{C}$ and $\mathcal{C}'$ around $0$, we have
\begin{align*}
a_n &= \frac{1}{2\pi i} \oint_{\mathcal{C}} \frac{A(\sqrt{z})}{z^{n+1}}\,dz = \frac{1}{2\pi i} \oint_{\mathcal{C}'}
 \sum_{h=0}^{\infty} \frac{t^h(1-t)^2}{(1-t^{h+1})^2} \cdot \frac{(1+t)^{2n+2}}{t^{n+1}} \cdot \frac{(1-t)\,dt}{(1+t)^3} \\
&=  \sum_{h=0}^{\infty} \frac{1}{2\pi i} \oint_{\mathcal{C}'}
\frac{(1-t)^3(1+t)^{2n-1}}{t^{n+1-h}(1-t^{h+1})^2}\,dt,
\end{align*}
using the substitution $\sqrt{z} = \frac{\sqrt{t}}{1+t}$ (or equivalently $z = \frac{t}{(1+t)^2}$). It follows that
\begin{align*}
a_n &= \sum_{h=0}^{\infty} [t^{n-h}] \frac{(1-t)^3(1+t)^{2n-1}}{(1-t^{h+1})^2} = [t^{n+1}] \sum_{h=0}^{\infty} \frac{t^{h+1}}{(1-t^{h+1})^2} (1-t)^3(1+t)^{2n-1} \\
&= [t^{n+1}] \sum_{k=1}^{\infty} \sigma(k) t^k (1-t)^3(1+t)^{2n-1} = \sum_{k=1}^{\infty} \sigma(k) [t^{n+1-k}] (1-t)^3(1+t)^{2n-1} \\
&= \sum_{k=1}^{n+1} \sigma(k) \Big( \binom{2n-1}{n+1-k} - 3 \binom{2n-1}{n-k} + 3 \binom{2n-1}{n-1-k} - \binom{2n-1}{n-2-k} \Big) \\
&= \sum_{k=1}^{n+1} \frac{4k \sigma(k) (2k^2-3n-2)(2n-1)!}{(n+1-k)!(n+1+k)!}\,.
\end{align*}
We remark here that the identity $[t^a] (1+t)^b = \binom{b}{a}$ that we are using even remains true for negative $a$ or for $a>b$ if the binomial coefficient is considered to be $0$ then. The manipulation in the final step is consistent with this.

Problem~\ref{prob:bona_dejonge} is equivalent to the inequality $C_{n+1}a_n > C_n a_{n+1}$ for $n \geq 3$, and since $C_{n+1} = \frac{4n+2}{n+2} C_n$, it can also be expressed as $(4n+2)a_n > (n+2)a_{n+1}$. Hence we are left to consider the inequality
\begin{equation*}
\sum_{k=1}^{n+1} \frac{4k \sigma(k) (2k^2-3n-2)(4n+2)(2n-1)!}{(n+1-k)!(n+1+k)!} > \sum_{k=1}^{n+2} \frac{4k \sigma(k) (2k^2-3n-5)(n+2)(2n+1)!}{(n+2-k)!(n+2+k)!}\,,
\end{equation*}
which reduces to
\begin{equation*}
\sum_{k=1}^{n+2} \frac{8k \sigma(k) (k^2-3n-4)(2k^2-n-2)(2n+1)(2n-1)!}{(n+2-k)!(n+2+k)!} < 0
\end{equation*}
after some manipulations. After multiplication by
\begin{equation*}
  n(n+1)(n+2)(2n+3) = \frac{(2n)(2n+2)(2n+3)(2n+4)}{8},
\end{equation*}
this can be expressed as
\begin{equation*}
\sum_{k=1}^{n+2} k \sigma(k) (k^2-3n-4)(2k^2-n-2) \binom{2n+4}{n+2-k} < 0.
\end{equation*}
Finally, replacing $n+2$ by $n$, what we have to prove in order to settle Problem~\ref{prob:bona_dejonge} is that
\begin{equation*}
\sum_{k=1}^n k \sigma(k) (k^2-3n+2)(2k^2-n) \binom{2n}{n-k} < 0
\end{equation*}
for all $n \geq 5$, which is precisely~\eqref{eq:main_ineq}.

\subsection{Asymptotic analysis}
\label{sec:asymptotics}

In the following, we provide the steps of the analysis of the sum $F(n)$, aided by the software package that was presented in the previous section. We will verify~\eqref{eq:main_ineq} for $n \geq N = 10000$ by means of an asymptotic analysis with explicit error terms. For $n < N$, one can verify the inequality with a computer by determining $F(n)$ explicitly in all cases.

All computations carried out in this section can be found in the SageMath
notebook located at
\begin{center}
      \url{https://arxiv.org/src/2403.09408v3/anc/2025-bona-dejonge.ipynb},
\end{center}
and a corresponding static version (containing computations and results)
is available at
\begin{center}
      \url{https://arxiv.org/src/2403.09408v3/anc/2025-bona-dejonge.html}.
\end{center}

\subsubsection*{Approximating the binomial coefficients}

It is useful to divide the entire sum by $\binom{2n}{n}$ and to approximate the quotient. Note that we have
\begin{equation*}
\frac{\binom{2n}{n-k}}{\binom{2n}{n}} = \frac{n(n-1) \cdots (n-k+1)}{(n+1)(n+2) \cdots (n+k)} = \frac{n}{n-k} \prod_{j=1}^k \frac{n-j}{n+j} = \frac{n}{n-k} \prod_{j=1}^k \frac{1-j/n}{1+j/n}.
\end{equation*}
This can be rewritten as
\begin{equation}\label{eq:binomial_exact}
\frac{\binom{2n}{n-k}}{\binom{2n}{n}} = \frac{n}{n-k} \exp \Big( \sum_{j=1}^k \log(1 - j/n) - \log(1+j/n) \Big) = \frac{n}{n-k} \exp \Big( {-} \sum_{j=1}^k \sum_{\substack{ r=1 \\ r \text{ odd}}}^{\infty} \frac{2j^r}{rn^r} \Big),
\end{equation}
an expression that will also be used later. It follows from it that
\begin{equation}\label{eq:binomial_crude}
\frac{\binom{2n}{n-k}}{\binom{2n}{n}} \leq \frac{n}{n-k} \exp \Big( {-} \sum_{j=1}^k \frac{2j}{n} \Big) \leq \frac{n}{n-k} \exp \Big( {-} \frac{k^2}{n} \Big).
\end{equation}
For small enough $k$, we can also obtain an asymptotic expansion in the same way as in the previous section. This will be discussed later.

\subsubsection*{The tails}

In order to replace the binomial coefficient by a simpler expression that is amenable to a Mellin analysis, we first have to handle the tails of the sum. 
For this purpose, we require an explicit bound for the divisor function $\sigma(k)$
in form of a constant $A > 0$ such that
\begin{equation}\label{eq:sigma-bound}
      \sigma(k) \leq A\cdot k\cdot \log\log n
\end{equation}
for $1\leq k\leq n$ when $n\geq N$. Assume temporarily that $N \leq k\leq n$.
Then, using an inequality due to Robin~\cite{Robin:1984:sigma-bounds},
\[
      \frac{\sigma(k)}{k} \leq e^{\gamma} \log\log k + \frac{0.6483}{\log\log k}
      \leq e^{\gamma} \log\log n + \frac{0.6483}{\log\log N}
      \leq \bigg(e^{\gamma} + \frac{0.6483}{\log\log N}\bigg) \log\log n.
\]
For $N = 10000$, we can choose $A = 52/25 \geq e^{\gamma} + 0.6483/\log\log N$,
and we can let the computer verify that~\eqref{eq:sigma-bound} also holds for
$1\leq k\leq n = N$.

For $k > \frac{n}{2}$,~\eqref{eq:binomial_crude} combined with the fact that the binomial coefficients $\binom{2n}{n-k}$ are decreasing in $k$ gives us
\begin{equation*}\frac{\binom{2n}{n-k}}{\binom{2n}{n}} \leq 2 e^{-n/4}.
\end{equation*}
So we have
\begin{align*}
0 &\leq \frac{1}{\binom{2n}{n}} \sum_{n/2 < k \leq n} k \sigma(k) (k^2-3n+2)(2k^2-n) \binom{2n}{n-k} \\
&\leq \sum_{n/2 < k \leq n} 4Ak^6 (\log \log n) e^{-n/4} \\
&\leq (A \log \log n) n^7 e^{-n/4},
\end{align*}
since it is easily verified that $\sum_{n/2 < k \leq n} k^6 \leq \frac{n^7}{4}$ for $n > 5$. Thus,
the contribution of the sum in this range is
\begin{equation}\label{eq:error-tails-greater-n2}
      \frac{1}{\binom{2n}{n}} \sum_{n/2 < k \leq n} k \sigma(k)(k^2 - 3n + 2)(2 k^2 - n)\binom{2n}{n-k} 
      = B_{n\geq N}\Big(\frac{52}{25} e^{-n/4} n^7 \log\log n\Big).
\end{equation}

Next, fix a constant $\alpha \in (\frac12,\frac34)$; the precise value is in principle irrelevant if one is only interested in an asymptotic formula (compare the choice of $\alpha$ in the previous section). However, for our computations with explicit error bounds it is advantageous to take a value close to $\frac34$, so we choose $\alpha = \frac7{10}$. We bound the sum over all $k \in [n^{\alpha},n/2]$. Here, we have
\begin{equation*}
\frac{\binom{2n}{n-k}}{\binom{2n}{n}} \leq 2 e^{-k^2/n}
\end{equation*}
by~\eqref{eq:binomial_crude}, thus (assuming that $N$ is large enough that $k^2 \geq n^{2\alpha} \geq 3n$ whenever $n \geq N$, which we can easily verify for $N = 10000$ and $\alpha = 7/10$)
\begin{align*}
0 &\leq \frac{1}{\binom{2n}{n}} \sum_{n^{\alpha} \leq k \leq n/2} k \sigma(k) (k^2-3n+2)(2k^2-n) \binom{2n}{n-k} \\
&\leq \sum_{n^{\alpha} \leq k \leq n/2} 4Ak^6 (\log \log k) e^{-k^2/n} \\
&\leq (4A \log \log n) \sum_{n^{\alpha} \leq k \leq n/2}  k^6 e^{-k^2/n}.
\end{align*}
The function $t \mapsto t^6 e^{-t^2/n}$ is decreasing for $t \geq \sqrt{3n}$, thus in particular for $t \geq n^{\alpha}$ under our assumptions. This implies that (by~\eqref{eq:integral-inequality})
\begin{equation*}
\sum_{n^{\alpha} \leq k \leq n/2}  k^6 e^{-k^2/n} \leq n^{6\alpha} e^{-n^{2\alpha-1}} + \int_{n^{\alpha}}^{\infty} t^6 e^{-t^2/n}\,dt.
\end{equation*}
The integral can be estimated by elementary means: for $T = n^{\alpha}$,
\begin{equation*}\int_{T}^{\infty} t^6 e^{-t^2/n}\,dt \leq \frac{1}{T} \int_{T}^{\infty} t^7 e^{-t^2/n}\,dt = \frac{n}{2T} \big(6 n^3+6 n^2 T^2+3 n T^4+T^6\big) e^{-T^2/n}.
\end{equation*}
For large enough $n \geq N$, this is negligibly small. This can be quantified 
with the help of some explicit computations with $B$-terms. We find that
\begin{equation}\label{eq:error-tails-nalpha-n2}
  \frac{1}{\binom{2n}{n}} \sum_{n^{\alpha} \leq k \leq n/2} \binom{2n}{n-k} k \sigma(k) (k^2 - 3n + 2) (2k^2 -n)
  = B_{n\geq N}\Big( \frac{25073}{5000} e^{-n^{2/5}} n^{\frac{9}{2}} \log\log n \Big).
\end{equation}

\subsubsection*{Approximating the summands}

So we are left with the sum over $k < n^{\alpha}$. Here, we can expand the exact expression in~\eqref{eq:binomial_exact}: this can be done
by cutting the sum over $r$ at some point (we choose the cutoff at $R = 9$) and estimating
\begin{align*} 
0 \leq \sum_{j=1}^k \sum_{\substack{ r=R \\ r \text{ odd}}}^{\infty} \frac{2j^r}{rn^r} &\leq \frac{2}{Rn^R(1-k^2/n^2)} \Big( k^R + \frac{k^{R+1}}{R+1} \Big) \\
&\overset{R=9}{=} B_{n\geq N}\Big( \Big(\frac{239}{10000} k^{10} + \frac{2223}{10000} k^9\Big) n^{-9} \Big),
\end{align*}
in the same way as in the previous section (compate~\eqref{eq:sumestimate}). This is followed by a Taylor expansion of the exponential multiplied with the expansion of $\frac{n}{n-k}$,
cf.~\eqref{eq:binomial_exact}.
The full and sufficiently precise asymptotic expansion can be found in our
auxiliary SageMath notebook. It reads
\begin{multline*}
\frac{\binom{2n}{n-k}}{\binom{2n}{n}} = 
e^{-k^2/n} \Big(
  1 - \frac{k^4 + k^2}{6 n^3} + \frac{k^8}{72 n^6}
  + \frac{k^2}{n^2} - \frac{k^{12}}{1296 n^9}
  - \frac{3 k^6}{20 n^5}
  + \frac{k^{16}}{31104 n^{12}} \\
  + \cdots + B_{n\geq N}\Big(\frac{k^{24}}{10000 n^{18}}\Big) 
  + B_{n\geq N}\Big(\frac{9 k^{21}}{10000 n^{16}}\Big) + \cdots
  \Big),
\end{multline*}
where the summands are ordered based on their individual upper growth bound
(found from substituting $k = n^{\alpha}$). The ellipses $\cdots$ indicate terms that are left out as the expression would otherwise be very long. If it were required,
this expansion could also be made more precise. Let us now split the
expression inside the brackets: let $S(n, k)$ denote the sum of all
``exact'' terms, and $S_B(n, k)$ the sum of all $B$-terms.
We want to evaluate
\begin{equation*}\sum_{1 \leq k < n^{\alpha}} k \sigma(k) \big( S(n,k) + S_B(n, k) \big) (k^2-3n+2)(2k^2-n) e^{-k^2/n}.
\end{equation*}
Let us first deal with the error estimate: since $|(k^2-3n+2)(2k^2-n)| \leq 2k^4 + 3n^2$ for all $k$ and $n$, it suffices to bound
\begin{align*}
&\sum_{1 \leq k < n^{\alpha}} k S_B(n, k) \sigma(k) (2k^4 + 3n^2) e^{-k^2/n} \\
&\hspace*{3cm}\leq (A \log \log n) \sum_{1 \leq k < n^{\alpha}} S_B(n, k) (2k^6 + 3n^2k^2) e^{-k^2/n} \\
&\hspace*{3cm}\leq (A \log \log n) \sum_{k \geq 1} S_B(n, k) (2k^6 + 3n^2k^2) e^{-k^2/n}.
\end{align*}
We can now apply~\eqref{eq:integral-inequality} to $t\mapsto t^{j} e^{-t^2/n}$ to find, with the help of
computer algebra,
\begin{equation}\label{eq:error-summation-approximation}
  \sum_{1\leq k < n^{\alpha}} k S_B(n, k) \sigma(k) (2k^4 + 3n^2) e^{-k^2/n} 
  = B_{n\geq N}\Big(\frac{146718899}{10000} \sqrt{n} \log\log n\Big).
\end{equation}
While this error is not quite as small as those collected so far, for $n = 10000$
it is still only about $26.1\%$ of the eventual main term.

Now we can consider the remaining sum
\begin{equation*}
\sum_{1 \leq k < n^{\alpha}} k \sigma(k) S(n,k)(k^2-3n+2)(2k^2-n) e^{-k^2/n}.
\end{equation*}
To this end, we first add back the terms with $k \geq n^{\alpha}$ and estimate their sum.

For $k\leq n^{3/4}$ the expansion in $S(n, k)$ can be bounded above,
$S(n, k) \leq c_1 \approx 4.372$, and for $k\geq n^{3/4}$
we have $S(n, k) \leq k^{20} / (10000 n^{15})$. As for the other factors in our summands,
we can bound $(k^2 - 3n + 2)(2k^2 -n)$ from above by $2k^4$. For an estimate of $\sigma(k)$
we use~\eqref{eq:sigma-bound} in the range $k < n^{3/4}$, and for the remaining case of
$k\geq n^{3/4}$ we use the well-known weaker bound $\sigma(k) \leq k^2$. This leaves us with
\begin{multline}\label{eq:error-completion-until-n34}
  \sum_{n^{\alpha} \leq k < n^{3/4}} k \sigma(k) S(n, k) (k^2 - 3n + 2)(2k^2 -n) e^{-k^2/n}\\ 
  \leq 2 A c_1 \log\log n \sum_{n^{\alpha} \leq k < n^{3/4}} k^6 e^{-k^2/n}
  = B_{n\geq N}\bigg( \frac{12553}{5000} e^{-n^{2/5}} n^{19/4} \log\log n \bigg),
\end{multline}
and
\begin{multline}\label{eq:error-completion-after-n34}
  \sum_{k \geq n^{3/4}} k \sigma(k) S(n, k) (k^2 - 3n + 2)(2k^2 -n) e^{-k^2/n}\\ 
  \qquad \leq \frac{4251}{5000\, n^{15}} \sum_{k \geq n^{3/4}} k^{27} e^{-k^2/n}
  = B_{n\geq N}\bigg( \frac{2867}{2500} e^{-n^{1/2}} n^{11/2} \bigg),
\end{multline}
where the sums have been bounded using the same integral estimate as before.

\subsubsection*{Mellin transform}

Having estimated all error terms related to pruning and completing the
tails of the sum, we now want to evaluate
\begin{equation}\label{eq:sum-for-mellin}
\sum_{k \geq 1} k \sigma(k) S(n,k) (k^2-3n+2)(2k^2-n) e^{-k^2/n}.
\end{equation}
This sum is a linear combination of sums of the form
\begin{equation}\label{eq:term-for-mellin}
\sum_{k \geq 1} k^a n^b \sigma(k) e^{-k^2/n}.
\end{equation}
In the precision chosen in our accompanying SageMath worksheet, there are 
121 such summands, to be precise. Set $t = n^{-1}$, refer to the sum in~\eqref{eq:term-for-mellin} as $g_{a, b}(t)$
and let $d_{a, b}$ denote the coefficients such that the sum in~\eqref{eq:sum-for-mellin} can be written 
as $\sum_{a, b} d_{a, b} g_{a, b}(t)$.
The Mellin transform of 
$g_{a, b}(t)$ is given by
\begin{equation*}g^*_{a, b}(s) = \int_0^{\infty} t^{s-1} \sum_{k \geq 1} k^a t^{-b} \sigma(k) e^{-k^2 t} \,dt = \zeta(2s-2b-a-1)\zeta(2s-2b-a) \Gamma(s-b).
\end{equation*}
By the Mellin inversion formula, the original function $g_{a,b}(t)$ can be recovered from its transform via
\begin{equation}\label{eq:bona-dejonge:inverse-mellin}
  g_{a, b}(t) = \frac{1}{2\pi i} \int_{c-i\infty}^{c +i\infty} \zeta(2s-2b-a-1)\zeta(2s-2b-a) \Gamma(s-b) t^{-s}\,ds
\end{equation}
for $c > \frac{a}{2} + b + 1$. We may shift the line of integration further left as long as we
collect all corresponding residues. In a first step, we shift the line of integration to
$c = 3/4$. While in some summands poles occur as far right as $s = 7/2$, 
a straightforward computation reveals that, as mentioned in the introduction of this article,
nontrivial cancellations take place: after summing all contributions, non-zero residues 
in the half-plane where $\Re(s) \geq 3/4$ can only be found for $s = 1$ and $s = 2$, where we
collect a contribution of
\begin{equation}\label{eq:asy-exact}
  \sum_{s_0\in \{1, 2\}} \sum_{a, b} d_{a, b}\res(g_{a, b}(s), s=s_0) = -\frac{1}{8t^2} + \frac{1}{24t} = -\frac{n^2}{8} + \frac{n}{24},
\end{equation}
which proves the asymptotic main term given in~\eqref{eq:asymptotics}.

To determine an explicit error bound for the shifted integrals, we can use
the exact same strategy as outlined in the example of the Ramanujan $Q$-function
that led us to the error bound in~\eqref{eq:ramanujan:shifted-integrals}.
Letting a computer work out the details, collect the contributions of
all summands, and combine them into a single $B$-term yields
\begin{multline}\label{eq:error-integral}
  \bigg\lvert \sum_{a, b} \frac{1}{2\pi i}\int_{3/4-i\infty}^{3/4+i\infty} g^{*}_{a, b}(s) t^{-s}~ds\bigg\rvert\\
  \qquad \leq\frac{1}{2\pi} \sum_{a, b} n^{c_{a, b}} \int_{-\infty}^{\infty} \lvert g^*_{a, b}(c_{a, b} + iw)\rvert~dw
  = B_{n\geq N}\Big(\frac{20320803}{5000} n^{3/4}\Big).
\end{multline}

\subsection{Collection of all contributions}

\begin{figure}[t]
  \begin{minipage}[t]{0.45\linewidth}
    \centering
    \includegraphics[width=\linewidth]{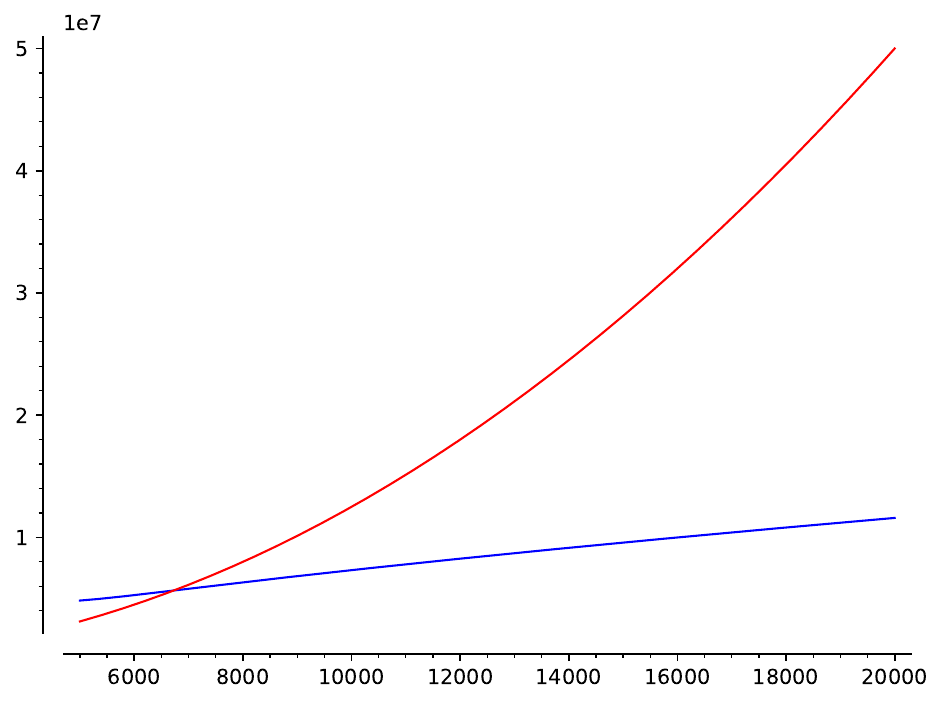}
  \end{minipage}\hfill
  \begin{minipage}[t]{0.45\linewidth}
    \centering
    \includegraphics[width=\linewidth]{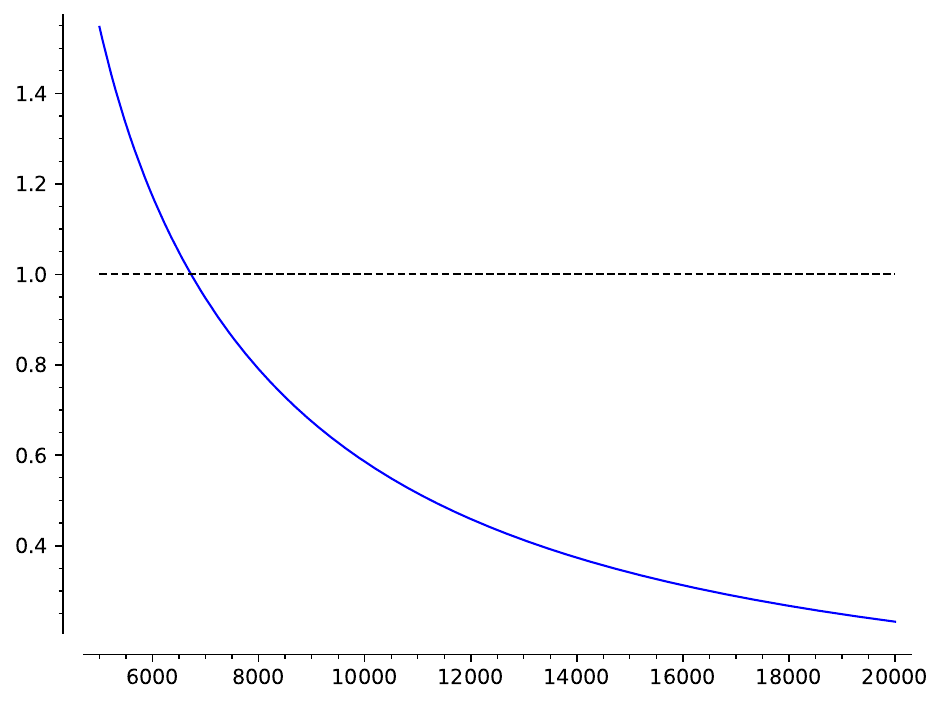}
  \end{minipage}
  \caption{
    Comparison of the absolute value of the asymptotic main term $-n^2/8 + n/24$ (red)
    against the accumulated total error (blue) on the left. The right plot depicts the
    ratio of the error bound to the main term.
  }\label{fig:error-comparison}
\end{figure}

Throughout Section~\ref{sec:asymptotics} we have accumulated several explicit error terms. They
are given in \eqref{eq:error-tails-greater-n2}, \eqref{eq:error-tails-nalpha-n2}, 
\eqref{eq:error-summation-approximation}, \eqref{eq:error-completion-until-n34},
\eqref{eq:error-completion-after-n34}, and \eqref{eq:error-integral}. Combining
them using crude estimates such as $\log\log n \leq n^{1/10}$ for $n\geq N$
proves the following theorem.

\begin{theorem}
  For $n\geq 10000$, the binomial sum $F(n)$ satisfies the asymptotic formula
  \[
      F(n) = \binom{2n}{n} \bigg(
            -\frac{n^2}{8} + \frac{n}{24}
            + B_{n\geq N}\Big(\frac{19374903}{2500} n^{3/4}\Big)
      \bigg).
  \]
\end{theorem}

Observe that for $n = 10000$ the certified error is already only
approximately $62.01\%$ of the absolute value of the exact main term.
Together with
the direct verification for $5\leq n < N$ this settles Problem~\ref{prob:bona_dejonge}. See
Figure~\ref{fig:error-comparison} for an illustration of the behavior of
the total error compared to the main term.

\subsection{Sketch of an alternative approach}

To conclude this paper, we briefly discuss an alternative approach that was kindly pointed out to us by a referee of the extended abstract~\cite{Hackl-Wagner:2024:binomial-mellin-AofA24}. Recall that the task is to prove the inequality~\eqref{eq:main_ineq}, i.e.,
\begin{equation*}
F(n) = \sum_{k=1}^n k \sigma(k) (k^2-3n+2)(2k^2-n) \binom{2n}{n-k} < 0.
\end{equation*}
Now one can use the well-known generating function identity 
\begin{equation*}
\sum_{n=k}^{\infty} \binom{2n}{n-k} x^n = x^k \sum_{m=0}^{\infty} \binom{2m+2k}{m} x^m = x^k \frac{1}{\sqrt{1-4x}} C(x)^{2k},
\end{equation*}
where $C(x) = \frac{1-\sqrt{4x}}{2x}$ is the generating function for the Catalan numbers, see e.g.~\cite[(5.72)]{Graham-Knuth-Patashnik:1994:concrete}. This gives an expression for $F(n)$ in terms of coefficients of functions involving $C(x)$. Specifically, we have 
\begin{align*}
F(n) &= 3n^2 \sum_{k=1}^n k \sigma(k) \binom{2n}{n-k} - n \sum_{k=1}^n k(7k^2+2) \sigma(k) \binom{2n}{n-k} \\
&\qquad + \sum_{k=1}^n 2 k^3 (k^2+2)  \sigma(k) \binom{2n}{n-k} \\
&= 3n^2 [x^n] \frac{1}{\sqrt{1-4x}} \sum_{k=1}^{\infty} k \sigma(k) (x C(x)^2)^k \\
&\qquad - n [x^n] \frac{1}{\sqrt{1-4x}} \sum_{k=1}^{\infty} k (7k^2+2) \sigma(k) (x C(x)^2)^k \\
&\qquad + [x^n] \frac{1}{\sqrt{1-4x}} \sum_{k=1}^{\infty} 2 k^3 (k^2+2)  \sigma(k) (x C(x)^2)^k \\
&= 3n^2 [x^n] \frac{1}{\sqrt{1-4x}} f_1(H(x)) - n [x^n] \frac{1}{\sqrt{1-4x}} f_2(H(x)) + [x^n] \frac{1}{\sqrt{1-4x}} f_3(H(x)),
\end{align*}
where $H(x) = xC(x)^2 = \frac{1-2x-\sqrt{1-4x}}{2x}$, and $f_1,f_2,f_3$ are given by the series
\begin{equation*}
f_1(z) = \sum_{k=1}^{\infty} k \sigma(k) z^k, \quad f_2(z) = \sum_{k=1}^{\infty} k (7k^2+2)  \sigma(k) z^k, \quad f_3(z) = \sum_{k=1}^{\infty} 2 k^3 (k^2+2) \sigma(k) z^k.
\end{equation*}
At the singularity $x = \frac14$, $H(x)$ has the expansion
\begin{equation*}
1 - 2 \sqrt{1-4x} + 2 (1-4x) + \cdots,
\end{equation*}
so we need the behavior of $f_1(z),f_2(z),f_3(z)$ around $z=1$. This can be determined by means of the Mellin transform: setting $z = e^{-t}$, we obtain for instance
\begin{equation*}
f_1(e^{-t}) = \sum_{k=1}^{\infty} k \sigma(k) e^{-kt},
\end{equation*}
whose Mellin transform is $\Gamma(s)\zeta(s-1)\zeta(s-2)$. Applying the inverse Mellin transform in the same way as in \eqref{eq:bona-dejonge:inverse-mellin} (though now with complex parameter $t$) yields
\begin{equation*}
f_1(e^{-t}) = \frac{\pi^2}{3t^3} - \frac{1}{2t^2} + O(t^K)
\end{equation*}
for any positive real $K$. This and analogous asymptotic formulas for $f_2$ and $f_3$ give us the behavior of $f_1(H(x))$,  $f_2(H(x))$ and $f_3(H(x))$ at the dominant singularity $\frac14$, from which the asymptotic formula~\eqref{eq:asymptotics} can be obtained by means of contour integration and singularity analysis. Carrying all of this out with explicit error terms comes with its own challenges, though, as one now has to deal with complex asymptotics.

\bibliographystyle{plainurl}
\bibliography{BonaDeJongeBTerms}

\end{document}